# Bivariate partial mapping for detecting causality in complex non-autonomous system

Yang Ni[a,#], Changqing Liu[a,#], Yifan Zhang[a], Yifan Gao[a], Haonan Guo[a], James Gao[a], Yingguang Li[a,*]

[a] College of Mechanical and Electrical Engineering, Nanjing University of Aeronautics and Astronautics, Nanjing, 210016, China.

[#] These authors contributed equally to this work.

[*] Corresponding author.

EMAIL: morganni@nuaa.edu.cn (Yang Ni), liuchangqing@nuaa.edu.cn (Changqing Liu), liyingguang@nuaa.edu.cn (Yingguang Li)

**Abstract**

Identifying causality is fundamental for human understanding of the world, where complex non-autonomous systems such as species population changes, brain activities, etc. are extensively existed. Since the phase spaces of such systems are not manifolds, the existing method based on convergent cross mapping is not applicable. This paper proposes a novel bivariate partial mapping method for detecting causality in complex non-autonomous systems. It transforms a non-autonomous system to an autonomous skew product system, and then, by considering the causality changes due to the transformation, detects causality of the original non-autonomous system from the transformed skew product system. The effectiveness of the proposed method is verified by mathematical cases and a real brain activity case, showing that the proposed method successfully detects the causality in complex non-autonomous systems.

## 1 Introduction

Complex non-autonomous systems are dynamical systems in which the variables are nonlinear coupled and the evolution rules are time-varying. Typical system cases include species population systems influenced by season changes, brain systems affected by external stimulus, etc., while detecting causality in such systems can be extremely difficult.

The coupling characteristic among variables of complex non-autonomous systems makes detecting methods based on Granger causality failure, because separability, the essential requirement of Granger causality, cannot be satisfied in such systems [1]. This characteristic also contradicts the underlying assumptions of Bayesian network models [2], making causal discovery methods based on which not applicable, such as methods by testing conditional independence [3], by testing graph matching properties [4] and by the structural asymmetry of causal functions [5–7].

Convergent cross mapping (CCM) provides a framework that detects causality from coupled variables based on nonlinear state space reconstruction [8]. The proposal of CCM represents a landmark achievement in detecting causality within complex systems. CCM reveals that if variables $x_i$ and $x_j$ in a system are causally linked, then they share a common attractor manifold $M$, and each variable is feasible to forecast the state of the other according to their shadow manifolds $M_{x_i}$ and $M_{x_j}$. The mathematical foundation for this implementation lies in that the shadow manifolds $M_{x_i}$ and $M_{x_j}$ are both diffeomorphic to the system manifold $M$, where $M_{x_i}$ and $M_{x_j}$ are reconstructed based on Takens Embedding Theorem [9]. Unfortunately, Takens Theorem is valid only when the system is autonomous [10]. For the non-autonomous system, system evolution trajectory may have intersections due to time-varying evolution rules. It means the phase space of the system may not be a manifold, as its intersection regions are not diffeomorphic to Euclidean space. Hence, Takens Theorem is not applicable for the non-autonomous system, leading to the failure of CCM. This problem can be solved by transforming the non-autonomous system to a skew product system, where the dynamics of time-varying evolution rules is included. In this way, phase space dimensions are expanded to exclude the intersections, making the phase space be a manifold, i.e., the system is autonomous [10].

However, it is easy to overlook that this transformation also changes system causality. The transformed skew product system expands phase space dimensions by introducing a new variable $\theta$, which is taken as the factor that drives the time-varying evolution rules of the original non-autonomous system. In fact, $\theta$ also disturbs the detecting of causality. The causality in the original non-autonomous system is generated only in $x_i$ and $x_j$, while that in the transformed skew product system is generated in $x_i$, $x_j$ and $\theta$. There are eight possible conditions according to whether $x_i$ is the cause of $x_j$, as well as the extent to which $\theta$ affects $x_i$ and $x_j$. Two conditions among these may be affected by $\theta$ if CCM is directly applied. (i) Let $x_i$ be the cause of $x_j$. $\theta$ has a strong effect on $x_i$ but a weak effect on $x_j$. Then it can hardly be taken that the reconstructed shadow manifold $M_{x_j}$ is diffeomorphic to the attractor manifold $M^{sp}$ of the skew product system. Therefore, the mapping of $x_i$ using $M_{x_j}$ will not converge for detecting true causality. (ii) Let $x_i$ not be the cause of $x_j$, but both $x_i$ and $x_j$ be strongly forced by $\theta$. Due to sufficient dynamics information of $\theta$ is contained, CCM tends to misidentify that $x_i$ is the cause of $x_j$, and this trend increases as the effect of $\theta$ becomes stronger. Actually, $\theta$ acts like a confounder when detecting causality in Condition (ii).

To this end, this paper proposes a novel bivariate partial mapping (BPM) method for detecting causality in complex non-autonomous systems. BPM transforms a non-autonomous system to an autonomous skew product system, and then, by considering the causality changes due to the transformation, detects causality of the original non-autonomous system from the transformed skew product system, as shown in Figure 1. Considering the conditions mentioned above, a bivariate shadow manifold $M_{x_j,\theta}$ is constructed, based on multivariate embedding theorem, to supplement the dynamics information of $\theta$ in $x_j$. Then, partial correlation analysis is applied between $x_i$ and the estimated $\hat{x}_i$ according to $M_{x_j,\theta}$, where $\theta$ is taken as the control variable to reduce its disturbance. Mathematical cases are used for verifications, showing that the proposed method can detect the causality of the original non-autonomous systems. A real brain activity case is also used for evaluation, where the causality detecting results are consistent with prior conclusions.

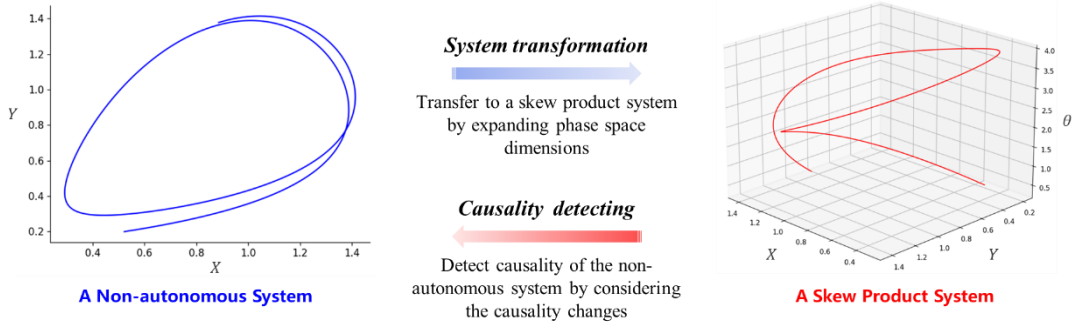

**Figure 1**. The overall idea of bivariate partial mapping method.

This paper is organized as follows. In Section 2, the methodology of this paper is provided. Section 3 shows the verifications while Section 4 presents discussions. Conclusions are given in Section 5.

**2 Methodology**

The methodology of BPM is introduced in this section, which is mainly divided into skew product system transformation and causality detecting.

*2.1 Skew product system transformation*

In the methods based on dynamical systems theory, delay embedding theorem serves as the mathematical foundation for causality detecting, which requires systems to be autonomous. Hence, the proposed BPM in this paper firstly transforms the complex non-autonomous system into a skew product system. The autonomy of the skew product system lies in it takes the external factor, driving the time-varying evolution rules of the original non-autonomous system, as one of the internal variables, which excludes the potential intersections of the system phase space by expanding the dimensions.

After expanding the dimensions of the original system phase space. A typical non-autonomous system can be represented as Equation (1):

$$\begin{cases} x_1(t+1) = f_1(t, x_1(t), \ldots, x_N(t)) \\ \ldots \\ x_N(t+1) = f_N(t, x_1(t), \ldots, x_N(t)) \end{cases} \quad (1)$$

where $x_1, \ldots, x_N$ are variables in the non-autonomous system, and $f_1, \ldots, f_N$ describe the corresponding evolution rules. It is time variable $t$ that makes the evolution rules of non-autonomous systems time-varying. That is, although in the same states, the subsequent evolution directions of the system may be various due to the differences of $t$. Therefore, the phase space of the non-autonomous system may have intersections and not be a manifold, leading to the failure of delay embedding theorem.

Let $t = \alpha\theta(t)$, $\alpha \neq 0$, the non-autonomous system can be transformed to a skew product system as Equation (2):

$$\begin{cases} x_1(t+1) = f_1(\alpha\theta(t), x_1(t), \ldots, x_N(t)) \\ \phantom{x_1(t+1) = f_1}\ldots \\ x_N(t+1) = f_N(\alpha\theta(t), x_1(t), \ldots, x_N(t)) \\ \phantom{xxxx}\theta(t+1) = \theta(t) + \frac{1}{\alpha} \end{cases} \quad (2)$$

It can be seen that the evolution of the skew product system only depends on system states, making the system autonomous. Comparing to the original non-autonomous system, variable $\theta$ is additionally included in the skew product system, which expands system phase space dimensions to exclude intersections.

## 2.2 Causality detecting

The skew product system transformation brings autonomy property for implementing delay embedding theorem, but the introduction of $\theta$ also changes system causality. As mentioned in Section 1, two of the eight conditions may lead to the false results of causality detecting if CCM is directly applied. The proposed BPM aims to detect the causality of the original non-autonomous system in the transformed skew product system, with the utilization of bivariate shadow manifold and partial correlation analysis.

### 2.2.1 Bivariate shadow manifold

Given that $x_j$ may lack sufficient $\theta$ dynamics information, BPM adopts bivariate shadow manifold, which is reconstructed using $x_j$ and $\theta$ based on Multivariate Embedding Theorem.

**Multivariate Embedding Theorem**. *Let $M$ be a compact $m$ dimensional manifold, $f: M \to M$ be a smooth flow, $H: M \to \mathbb{R}^n$ be a complete smooth observation function. Then it is a generic property that for $\forall x \in M$, delay mapping*

$$\Psi(x) = \left(H(x), H(f(x)), \ldots, H(f^{d-1}(x))\right) \in \mathbb{R}^{nd} \quad (3)$$

*is an embedding when $nd \geq 2m + 1$.*

The detailed proof of Multivariate Embedding Theorem can be seen in S1. Note that when $n = 1$, Multivariate Embedding Theorem is equivalent to Takens Theorem [9].

According to Multivariate Embedding Theorem, the reconstructed bivariate shadow manifold $M_{x_j,\theta}$ is generically diffeomorphic to the skew product system manifold $M^{sp}$, which is used for causality detecting.

*2.2.2 Partial correlation analysis*

Given that any two variables $x_i$ and $x_j$ may be strongly affected by $\theta$ that acts as a confounder, analysis based on Pearson correlation coefficient between $x_i$ and the estimated $\hat{x}_i$ according to $x_j$ is likely to misidentify the causality. Actually, as the effects of $\theta$ increase, the intrinsic dynamics of $x_i$, $x_j$ gradually become subordinate to $\theta$, i.e., the evolution of $x_i$, $x_j$ and $\theta$ are gradually linearly related [11]. It raises the possibility of the misidentification.

Therefore, BPM takes $\theta$ as control variable, and makes partial correlation analysis of $x_i$ and $\hat{x}_i$ as follows:

$$Pcc(x_i, \hat{x}_i|\theta) \tag{4}$$

where $Pcc$ is partial correlation coefficient. Theoretically, partial correlation analysis can exclude the linear effect of $\theta$ in $x_i$ and $\hat{x}_i$, which reduces the effects of $\theta$ for detecting the causality of the original non-autonomous system.

*2.2.3 An algorithm for bivariate partial mapping*

The algorithm of BPM is similar to that of CCM, where simplex projection [12] is applied, involving exponentially weighted distances from nearby points on a reconstructed manifold to do kernel density estimation [9].

For $\theta$ and any two variables $x_i, x_j$ within a transformed skew product system, their temporal data of length $L$ can be collected as $\{\theta\} = \{\theta(1), \theta(2), \dots, \theta(L)\}$, $\{x_i\} = \{x_i(1), x_i(2), \dots, x_i(L)\}$ and $\{x_j\} = \{x_j(1), x_j(2), \dots, x_j(L)\}$. Based on Multivariate Embedding Theorem, shadow manifold $M_{x_i}$ and bivariate shadow manifold $M_{x_j,\theta}$ can be reconstructed, which are constituted with the sets of lagged-coordinate vectors

$$\boldsymbol{x}_i(t) = <x_i(t), x_i(t-\tau), \dots, x_i(t-(E-1)\tau)> \tag{5}$$

and

$$\boldsymbol{x_j, \theta}(t) = <x_j(t), \theta(t), \dots, x_j\left(t-\left(\frac{E}{2}-1\right)\tau\right), \theta\left(t-\left(\frac{E}{2}-1\right)\tau\right)> \tag{6}$$

$\tau$ denotes the positive time lag of smooth flow $f$. $E$ is the dimension of embedded state space. $t$ is the integer value from $1 + (E-1)\tau$ to $L$.

$E + 1$ nearest neighbors of $x_j, \theta(t)$ can be found by $t_1, t_2, \ldots, t_{E+1}$ on $M_{x_j, \theta}$.

These time indices corresponding to nearest neighbors to $x_i(t)$ on $M_{x_i}$ are then used for estimating $\hat{x}_i$ from $E + 1$ $x_i(t_k)$ values as follows:

$$\hat{x}_i^{M_{x_j,\theta}} = \sum_{k=1}^{E+1} w_k \, x_i(t_k) \tag{7}$$

where $k = 1, 2, \ldots, E + 1$. $w_k$ is a weight calculated according to the distance between $x_j, \theta(t)$ and its $k^{th}$ nearest neighbor on $M_{x_j, \theta}$.

$$w_k = \frac{u_k}{\sum_{K=1}^{E+1} u_K} \text{ w.r.t } u_k = e^{-\frac{d\{[x_j,\theta(t)],[x_j,\theta(t_k)]\}}{d\{[x_j,\theta(t)],[x_j,\theta(t_1)]\}}} \tag{8}$$

where $d$ denotes the measurement of Euclidean distance.

If $x_i$ is the cause of $x_j$, $\hat{x}_i$ should converge to $x_i$ under the control of $\theta$ as $L$ increases, and BPM adopts partial correlation coefficient calculated in Expression (4) as the index for detecting the causality of the original non-autonomous system.

## 3 Verifications

Verifications are carried out on mathematical cases and a real brain activity case, where the classic CCM in literature is taken for comparing. The convergence of the correlation coefficient in BPM and CCM is used to distinguish the causality.

### *3.1 Verifications on mathematical cases*

Lotka-Volterra model that shows the changes in species population is selected for verifying, where time variable $t$ is explicitly introduced to reflect the disturbance caused by external factors, such as food quantity, temperature, humidity, etc. As a result, the modified Lotka-Volterra models belong to complex non-autonomous systems.

Three cases of the modified Lotka-Volterra models are provided to show the effectiveness of BPM. Expression (9)-(14) give the mathematical expressions of the non-autonomous system and corresponding skew product system of each case, where the skew product system is transformed by letting $t = \frac{4}{3}\theta[t]$.

**Case 1**

The dynamical equations of the non-autonomous system:

$$\begin{cases} X[t+1] = X[t](3.8(1 - X[t]) - 0.08(1 + \beta_1 \cos\left(\frac{3}{2}\pi t\right))Y[t]) \\ Y[t+1] = Y[t](3.5(1 - Y[t]) - 0.08(1 + \beta_2 \cos\left(\frac{3}{2}\pi t\right))X[t]) \end{cases} \tag{9}$$

The dynamical equations of the skew product system:

$$\begin{cases} X[t+1] = X[t](3.8(1-X[t]) - 0.08(1+\beta_1\cos(2\pi\theta[t]))Y[t]) \\ Y[t+1] = Y[t](3.5(1-Y[t]) - 0.08(1+\beta_2\cos(2\pi\theta[t]))X[t]) \\ \theta[t+1] = \theta[t] + \frac{3}{4} \text{ (mod 1)} \end{cases} \quad (10)$$

**Case 2**

The dynamical equations of the non-autonomous system:
$$\begin{cases} X[t+1] = X[t](3.6 + \beta_3\cos(\frac{3}{2}\pi t))(1-X[t]) \\ Y[t+1] = Y[t]((3.5 + \beta_4\sin(\frac{3}{2}\pi t))(1-Y[t]) - 0.08X[t]) \end{cases} \quad (11)$$

The dynamical equations of the skew product system:
$$\begin{cases} X[t+1] = X[t](3.6 + \beta_3\cos(2\pi\theta[t]))(1-X[t]) \\ Y[t+1] = Y[t]((3.5 + \beta_4\sin(2\pi\theta[t]))(1-Y[t]) - 0.08X[t]) \\ \theta[t+1] = \theta[t] + \frac{3}{4} \text{ (mod 1)} \end{cases} \quad (12)$$

**Case 3**

The dynamical equations of the non-autonomous system:
$$\begin{cases} X[t+1] = X[t](3.8 + \beta_5\sin(\frac{3}{2}\pi t))(1-X[t]) \\ Y[t+1] = Y[t](3.5 + \beta_6\sin(\frac{3}{2}\pi t))(1-Y[t]) \end{cases} \quad (13)$$

The dynamical equations of the skew product system:
$$\begin{cases} X[t+1] = X[t](3.8 + \beta_5\sin(2\pi\theta[t]))(1-X[t]) \\ Y[t+1] = Y[t](3.5 + \beta_6\sin(2\pi\theta[t]))(1-Y[t]) \\ \theta[t+1] = \theta[t] + \frac{3}{4} \text{ (mod 1)} \end{cases} \quad (14)$$

The causality of Case 1-3 can be described as in Figure 2. Note that in the above cases, $\theta$ is included in "$\cos(2\pi\theta[t])$" or "$\sin(2\pi\theta[t])$", making the influence of $\theta$ periodic. Hence, the dynamical equation of $\theta$ can take "mod 1" for simplification.

Figure 3 (A) shows the detecting results of CCM and BPM for Case 1, where $X$ and $Y$ has a bidirectional causality. When $\beta_1 = 1, \beta_2 = 0$, detecting whether $Y$ is the cause of $X$ falls under Condition (i). It can be seen that CCM incorrectly identifies $Y$ is not the cause of $X$, as the blue line converges to around 0. It is successfully detected by BPM as the orange line converges to over 0.6. Besides, just like CCM, BPM also correctly detects that $X$ causes $Y$, as the purple and red lines converge to nearly 1.0.

Figure 3 (B) shows the detecting results of CCM and BPM for Case 2, where $X$ unidirectionally causes $Y$. When $\beta_3 = 0.3, \beta_4 = 0.5$, detecting whether $Y$ is the cause of $X$ belongs to Condition (ii). CCM incorrectly identifies $Y$ is the cause of $X$, as the blue line converges nearly to 1. On the contrary, the orange line converges to nearly 1, showing that BPM makes the right causality identification. Both BPM and CCM identify that $X$ is the cause of $Y$, as the purple and red lines converge to nearly 1.

Figure 3 (C) shows the detecting results of CCM and BPM for Case 3, where no causality exists between $X$ and $Y$. When $\beta_5 = 0.2, \beta_6 = 0.2$, detecting whether $Y$ is the cause of $X$ belongs to Condition (ii). Likewise, CCM incorrectly identifies $Y$ causes $X$ as the blue line converges to nearly 1, while BPM makes the right causality identification as the purple and orange lines converge to around 0.

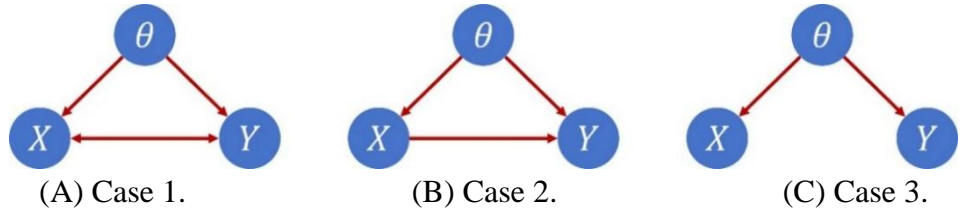

(A) Case 1.  (B) Case 2.  (C) Case 3.

**Figure 2**. The causality of Case 1-3.

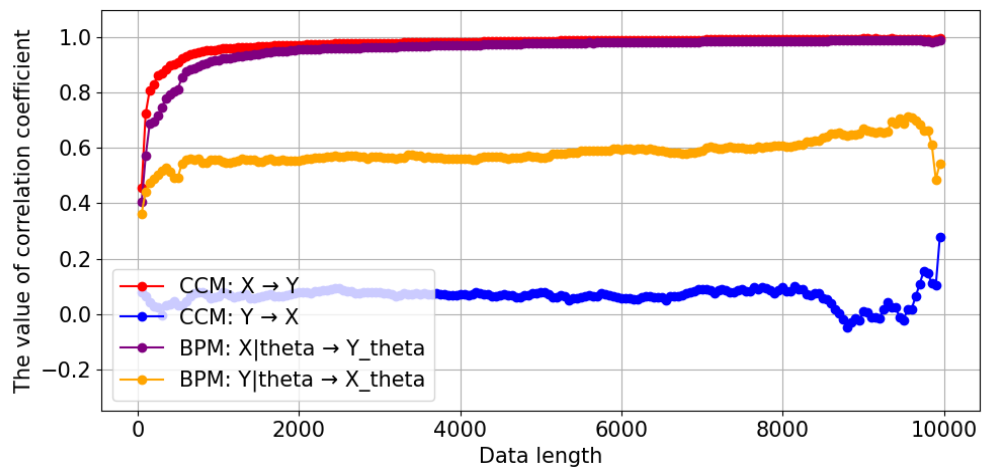

(A) Case 1.

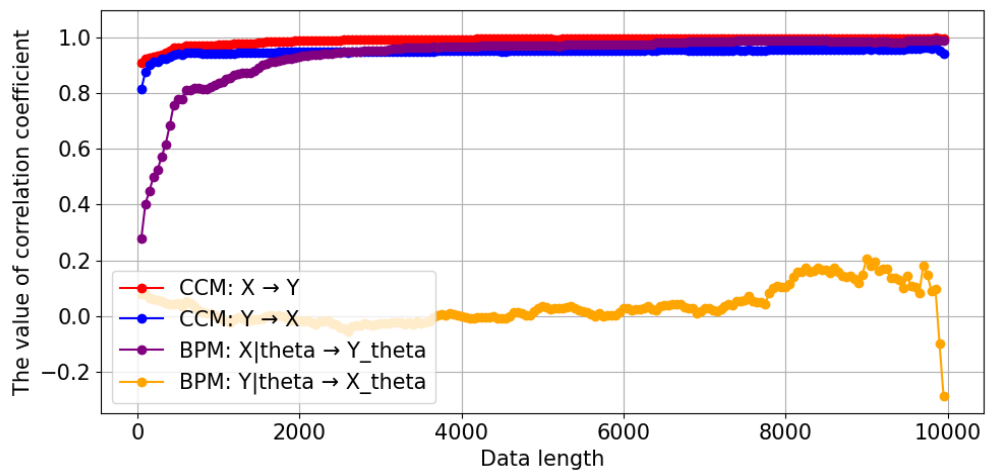

(B) Case 2.

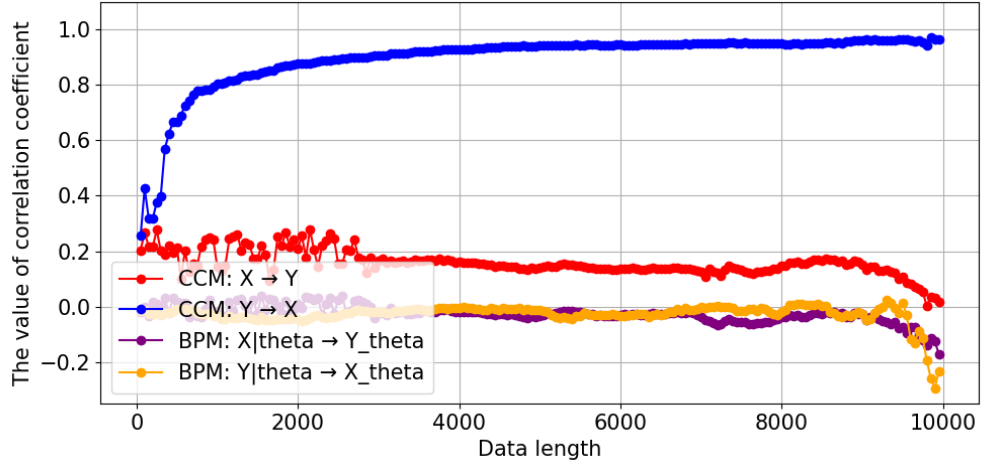

(C) Case 3.

**Figure 3**. The detecting results of CCM and BPM for mathematical Case 1-3.

*3.2 Verifications on a brain activity case*

The proposed BPM is also performed on a brain activity case, which investigates the causality of the functions of two brain regions in a mouse involving visual stimulus. In a trial, the mouse responds to the visual stimulus and receives corresponding rewards, where Neuropixels probes are inserted to collect brain signals. Over multiple trials, a 1.5-h long firing rate time series is recorded to reflect the activity degree of brain regions. Details are provided in Literature [13], and data are available at https://data.internationalbrainlab.org.

*3.2.1 Preliminary*

The caudate putamen (CP) and primary somatosensory area (SSp) of mice are selected as the objects of this case, which are causally linked according to the research conclusion in Literature [14]. The rewards given to the mouse can help it gradually establish a conditioned reflex, affecting the evolution rules of brain signals. Hence, the system composed of CP and SSp signals can be regarded as a complex non-autonomous system.

The sampling frequency of the probes is set to be 80-Hz. The duration of a single trial is between 3 and 4 seconds, and the 1.5-h long time series is divided in nonoverlapping 10-s segments. According to Literature [13], each segment is taken as having the same level of conditioned reflex.

*3.2.2 Causality detecting*

The non-autonomous system composed of CP and SSp signals is firstly transformed to a skew product system by introducing variable $\theta$, which indicates the level of conditioned reflex that the mouse is at. According to the preliminary setting that each segment has the same level of conditioned reflex, $\theta$ can be set to have the same value in each segment and varies across segments for simplification as follows:

$$\begin{cases} \theta[t+1] = \theta[t] & if\ t+1\ and\ t\ are\ within\ the\ same\ segment \\ \theta[t+1] = \theta[t] + \delta & if\ t+1\ and\ t\ are\ not\ in\ the\ same\ segment \end{cases} \quad (15)$$

where $\delta$ is a constant and set to be 1 in this paper.

For each trial, the brain signals are evoked by the visual stimulus at the starting moment, which makes the signals discontinuous. To satisfy the requirements of delay embedding theorem that the flow and observation function need to be smooth, the temporal data for making delay mapping are sampled from the same trial. This does not affect the reconstruction of shadow manifolds, as the delay mapping using data from different trials can be regarded as the embedding carried out on the different positions of the skew product system manifold.

Figure 4 shows the causality detecting results using CCM and BPM respectively. It can be seen that the detecting results of CCM does not clarify the causality between CP and SSp, as the purple and orange lines converge to around 0.3 and 0.1 respectively. On the contrary, BPM indicates that there is a bidirectional causality between CP and SSp considering the purple and orange lines converge to around 0.5 and 0.7 respectively. The detecting results of BPM are consistent with the prior conclusion, demonstrating the effectiveness of the proposed method.

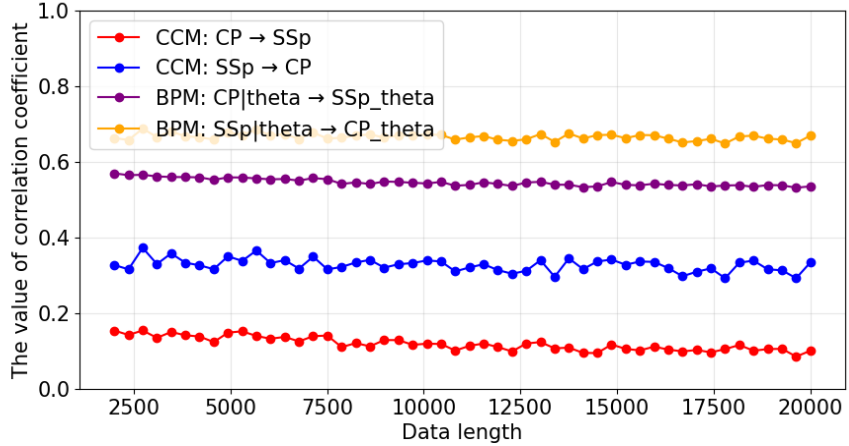

**Figure 4**. The detecting results of CCM and BPM for the brain activity case.

## 4 Discussions

Discussions mainly proceed from the following two aspects.

In the aspect of causality detecting, the main reason of the failure of CCM lies in the detecting results belong to the transformed skew product system rather than the original non-autonomous system. This paper does not aim to raise any new doubts about the validity of CCM, but emphasizes that the introduced variable $\theta$ can change system causality. Hence, by considering the causality changes due to system transformation, the proposed BPM is designed to detect the causality of the original non-autonomous system in the transformed skew product system.

In the aspect of $\theta$, its dynamics can be determined using $\theta(t+1) = \theta(t) + \frac{1}{\alpha}$, and it can be further simplified according to the prior knowledge of systems, such as the periodicity prior in the mathematical cases and the conditioned reflex level prior in the brain activity case of this paper. Note that these simplifications can introduce the prior knowledge into the reconstruction of bivariate shadow manifold and the partial correlation analysis, which improves the accuracy of BPM. Besides, the value of $\theta$ should theoretically have upper and lower bounds to keep system manifolds compact, satisfying the requirements of delay embedding theorem. Otherwise, the theorem might fail, leading to incorrect results of causality detecting.

## 5 Conclusions

This paper proposed BPM for detecting causality in complex non-autonomous systems. The mathematical cases shows that BPM can be applied to more conditions than CCM, and a brain activity case also makes further verifications.

More application scenarios of BPM will be investigated in the future, including brain map reconstructions, etiology analyses, manufacturing process modeling, etc.

**Acknowledgement**

The reported research was funded by the National Natural Science Foundation of China (grant No. 52505561, grant No. 52575579).

**S1 The proof of Multivariate Embedding Theorem**

According to Whitney Embedding Theorem, let $M$ be a m-dimensional smooth compact manifold, then it can be smoothly embedded into $\mathbb{R}^{2m+1}$, hence $nd$ is required to be greater than or equal to $2m + 1$. Then, to show that $\Psi: M \to \mathbb{R}^{nd}$ is an embedding, it is to show that $\Psi$ is an immersion as well as a homeomorphic mapping.

Proof of immersion

$\Psi$ is an immersion is equivalent to that the Jacobi matrix of $\Psi$ that denotes as $D\Psi$ is of full rank:

$$\text{rank}(D\Psi) = dimM \tag{16}$$

where $\Psi(x) = \left(H(x), H(f(x)), \ldots, H(f^{d-1}(x))\right)$. Observation function $H: M \to \mathbb{R}^n$ is complete, which means any $x \in M$ can be uniquely determined according to its delayed observation, and $H$ can distinguish system trajectories. Therefore, for any two points $x_a, x_b$ of $M$, it can be deduced that $H(f^l(x_a)) = H(f^l(x_b))$ if and only if $x_a = x_b$, where $l = 0, 1, 2, \ldots, d-1$. The following part shows that if observation function $H$ is incomplete, it will contradict with this deduction.

Let $T_x$ denote the point $x$ within a tangent space. $DH: T_xM \to T_x\mathbb{R}^n$ shows the ability of information retention by $H$ at the local region of point $x$. Then $rank(DH) < dimM$ indicates that when $x$ changes in a specific dimension on $M$, the observed value of $H$ remains unchanged, as the dynamics of this dimension is lacked in $H$.

Assuming that there is a point $x_0$ letting $rank(DH) < dimM$, then for point $u$ in the neighborhood $U \subset M$ of point $x_0$, value $H(f^l(u))$ is only determined by $k < m$ dimensions within the coordinate vector $<u_1, ..., u_m>$. When the values of the remaining $m - k$ dimensions are different, $H(f^l(u_a)) = H(f^l(u_b))$ will lead to $u_a \neq u_b$, which is inconsistent to the above deduction.

Therefore, the completeness of observation function can derive that $rank(DH) = dimM$ and, further, $rank(D\Psi) = dimM$. $\Psi$ is an immersion.

Proof of homeomorphic mapping

$\Psi$ is a homeomorphic mapping if and only if $\Psi$ is bijective and its forward/inverse mapping is continuous.

The bijectivity of $\Psi$ requires it to be injective and surjective. For injectivity, it can be proved according to the deduction that $H(f^l(x_a)) = H(f^l(x_b))$ derives $x_a = x_b$. For surjectivity, it can be proved as any point in the image set of $\Psi(x)$ has a corresponding point in $M$. Therefore, $\Psi$ is bijective.

The continuity of $\Psi$ can be deduced considering that flow $f$ and observation function $H$ within $\Psi$ are both smooth. For its inverse, it can be proved as follows:

1. As $M$ is compact, any open subset $C \subset M$ is compact.

2. Due to $\Psi$ is continuous. $\Psi(C) \subset \mathbb{R}^{nd}$ is compact.

3. $\Psi(C)$ locates in Hausdorff space. Given $C$ is open and $\Psi(C)$ is compact, $\Psi(C)$ is also a open set.

4. According to $\Psi$ is bijective, its inverse $\Psi^{-1}$ maps open set $\Psi(C) \subset \mathbb{R}^{nd}$ to open set $C \subset M$. It can deduce that $\Psi^{-1}$ is continuous.

Therefore, $\Psi$ is bijective and its forward/inverse mapping is continuous, and it is a homeomorphic mapping.

Conclusion

$\Psi$ is an immersion as well as a homeomorphic mapping, showing that $\Psi: M \to \mathbb{R}^{nd}$ is an embedding.

**References**

[1]   C.W.J. Granger, Investigating Causal Relations by Econometric Models and Cross-spectral Methods, Econometrica 37 (1969) 424–438.

[2]   J. Pearl, Causality, Cambridge University Press, New York, 2009.

[3]   P. Spirtes, C. Glymour, R. Scheines, Causation, Prediction, and Search, The MIT Press, New York, 1993. https://doi.org/10.1007/978-1-4612-2748-9.

[4]   D.M. Chickering, Optimal Structure Identification With Greedy Search, Journal of Machine Learning Research 3 (2002) 507–554.

[5]   S. Shimizu, P.O. Hoyer, A. Hyvärinen, A linear non-gaussian acyclic model for causal discovery, Journal of Machine Learning Research 7 (2006) 2003–2030.

[6]   P. Daniušis, D. Janzing, J. Mooij, J. Zscheischler, B. Steudel, K. Zhang, B. Schölkopf, Inferring deterministic causal relations, in: UAI, 2010.

[7]   R. Cai, J. Qiao, K. Zhang, Z. Zhang, Z. Hao, Causal Discovery with Cascade Nonlinear Additive Noise Models, in: IJCAI, 2019: pp. 1609–1615.

[8]   G. Sugihara, R.M. May, H. Ye, C.-H. Hsieh, E. Deyle, M. Fogarty, S. Munch, Detecting Causality in Complex Ecosystems, Science (1979). 338 (2012). https://www.science.org.

[9]   F. Takens, Dynamical Systems and Turbulence, Springer, New York, 1981.

[10]  J. Stark, Delay Embeddings for Forced Systems. I. Deterministic Forcing, J. Nonlinear Sci 9 (1999) 255–332.

[11]  K. Josi, Synchronization of chaotic systems and invariant manifolds, Nonlinearity 13 (2000) 1321–1336.

[12]  G. Sugihara, R.M. May, Nonlinear forecasting as a way of distinguishing chaos from measurement error in time series, 344 (1990) 734–741.

[13]  International Brain Laboratory, et al, A brain-wide map of neural activity during complex behaviour, Nature 645 (2025) 177–191. https://doi.org/10.1038/s41586-025-09235-0.

[14]  N.A. Steinmetz, P. Zatka-Haas, M. Carandini, K.D. Harris, Distributed coding of choice, action and engagement across the mouse brain, Nature 576 (2019) 266–273. https://doi.org/10.1038/s41586-019-1787-x.